\documentclass[11pt,a4paper]{article}

\usepackage[T1]{fontenc}
\usepackage[utf8]{inputenc}
\usepackage{lmodern}

\usepackage{amsmath,amssymb,amsfonts,mathtools}
\usepackage{bm}

\usepackage{amsthm}

\theoremstyle{definition}

\usepackage{graphicx}
\usepackage{booktabs}
\usepackage{caption}
\usepackage{subcaption}

\usepackage{siunitx}

\usepackage[
a4paper,
left=3cm,
right=3cm,
top=2.5cm,
bottom=2.5cm
]{geometry}

\usepackage{authblk}

\usepackage[
    backend=biber,
    style=numeric-comp,
    sorting=none,
    maxbibnames=99
]{biblatex}

\usepackage{hyperref}

\title{%
    A Comparison of Active Flux Methods for the Vlasov-Poisson System
}

\author[1]{L. Hensel}
\author[3]{Y. Kiechle}
\author[2]{R. Grauer}
\author[2]{G. Gr\"unwald}
\author[3]{C. Helzel}
\author[1]{K. Kormann}

\affil[1]{Numerical Mathematics, Ruhr-Universit\"at Bochum, Universit\"atsstra\ss{}e 150, D-44801 Bochum, Germany}
\affil[2]{Theoretical Physics I, Ruhr-Universit\"at Bochum, Universit\"atsstra\ss{}e 150, D-44801 Bochum, Germany}
\affil[3]{Faculty of Mathematics and Natural Sciences, Institute of Mathematics, Heinrich-Heine-University D\"usseldorf, Germany}

\date{}

\begin{document}

\maketitle

\begin{abstract}
	Active Flux is a third-order accurate, fairly novel finite volume method for hyperbolic conservation laws that is becoming increasingly popular.
	It evolves additional nodal degrees of freedom (DOF) located on cell interfaces and shared by neighboring cells. The numerical fluxes are then computed from these DOFs.
	A crucial component of Active Flux methods is the evolution operator of the point values, which enables the natural use of semi-Lagrangian ideas and makes Active Flux an attractive candidate for a grid-based approach to the Vlasov equation.
	Here, we compare two recently proposed Active Flux methods for the 1D1V Vlasov--Poisson system: a split-step method and an unsplit method.
\end{abstract}

\vspace{-0.5cm}
\section{Introduction}
\label{sec:1}

The Vlasov--Poisson system describes the time evolution of the particle distribution function $f(x,v,t)$ in phase space for an electrostatic collisionless plasma. It couples a linear transport equation in phase space, to the Poisson equation, from which the self-consistent electric field is obtained. The source for the Poisson equation is given by the plasma's charge density, which is computed as the zeroth velocity moment of the distribution function. The equation for the electron distribution function in a static neutralizing ion background reads:
\begin{equation}
	\begin{aligned}
		\partial_t f(x,v,t) + &v \partial_x f(x,v,t) - E(x,t) \partial_v f(x,v,t) = 0 \\
		-\Delta \phi(x,t) &= 1 - \int_{\mathbb{R}} f(x,v,t)\,dv, \quad E = - \nabla \phi. 
	\end{aligned}
\end{equation}
It represents the one-dimensional electrostatic approximation of the Vlasov--Maxwell equations.
Numerical modelling of both Vlasov--Maxwell and Vlasov--Poisson is still an active field of research in an astrophysical \cite{allmann2022energy} as well as fusion context. The main methods are Particle in Cell (PIC) methods \cite{Kraus_Kormann_Morrison_Sonnendrucker_2017}, grid-based methods such as discontinuous Galerkin \cite{ROSSMANITH20116203}, finite volume methods \cite{BanksHittinger2010} and semi-Lagrangian schemes \cite{FSB2001, crouseilles_conservative_2010}. A significant proportion of these approaches are founded upon operator or dimensional splittings, which were first established in Cheng and Knorr's work \cite{cheng1976integration}. 
Later, different splitting techniques have been developed to achieve higher order, preserve physical quantities \cite{YOSHIDA1990262, CHENG2014145} or to maximize efficiency \cite{blanes2002practical}. 
In \cite{mclachlan2021tuning} optimized splitting schemes (with higher cost) showed strong efficiency advantages over more popular schemes with minimal number of steps.
Active Flux methods are finite volume methods that enjoy increased popularity in the last years. In 2011, Roe and his former students began working on this method and highlighted in a number of publications the advantages of Active Flux methods \cite{eymann2011activeflux, fan2015investigations, roe2018comparing}. Specifically, the method is third-order accurate and has a compact stencil in both space and time. These features enable the method to be highly parallelizable and achieve high accuracy on coarse grids due to its very low dissipation. 
Initially, Active Flux was proposed as a fully discrete method, an approach that is also considered in  \cite{chudzik2025fullydiscretetrulymultidimensional}.
In recent publications, Active Flux methods for the Vlasov--Poisson \cite{HENSEL2025114294, KIECHLE2025113693, procKCH2023} and the Vlasov--Maxwell system \cite{grunwald2025solving} have been developed. The objective of this article is to provide a comparative study of two approaches, namely a split-step and an unsplit method, for simulating the 1D1V Vlasov--Poisson system.
The article is structured as follows: Section \ref{sec:2} reviews the two methods followed by a theoretical comparison. A comparative numerical study is given in Section \ref{sec:3}.

\section{Active Flux based Vlasov Solvers}
\label{sec:2}

\subsection{General Structure of Active Flux Schemes on Cartesian Grids}
Active Flux adopts a globally continuous reconstruction by introducing point value degrees of freedom at cell interfaces that are shared between neighboring cells. This eliminates the need for Riemann solvers and leads to compact stencils in both space and time, closely following Roe's principle that a CFD algorithm should "use all of the information on which the solution depends, and should use only that information" \cite{Roe_Musings_2024}. 
Therefore, the aim is to update the DOFs using information from a domain as close as possible to the physical domain of dependence. A key distinction from classical finite volume methods lies in the computation of fluxes: the additional interface degrees of freedom are used within high-order quadrature rules to evaluate fluxes across cell boundaries directly from the continuous reconstruction. While the evolution of the cell averages remains conservative, the update of the interface point values is not restricted to conservative formulations. Instead, it can rely on update operators derived from the exact solution of the underlying initial value problem.
For linear advection, this naturally results in a semi-Lagrangian or characteristic tracing approach, in which the continuous, piecewise quadratic reconstruction is evaluated at the foot point of the characteristic.

\subsection{An unsplit Active Flux Method for Vlasov--Poisson}
The following subsection provides a brief overview of the unsplit Active Flux method described in \cite{KIECHLE2025113693}. 

\begin{figure}
	\centering
	\begin{subfigure}[b]{0.49\textwidth}
		\centering
		\includegraphics[width=\textwidth]{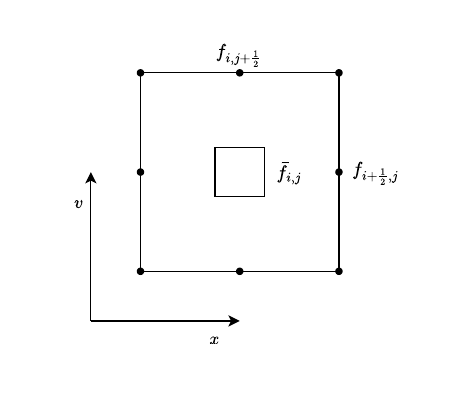}
		\caption{Cell for the unsplit method with point values on the interfaces and a two-dimensional cell average}
		\label{fig:y equals x}
	\end{subfigure}
	\hfill
	\begin{subfigure}[b]{0.49\textwidth}
		\centering
		\includegraphics[width=\textwidth]{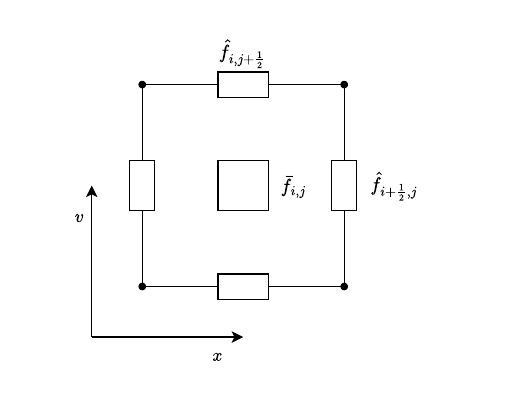}
		\caption{Cell for the split-step method with nodal-point values, 1D line averages over the edges and a two-dimensional cell average}
		\label{fig:three sin x}
	\end{subfigure}
	\hfill
	\caption{Different grid structures for the unsplit and split-step Active Flux methods}
	\label{fig:grid_cells}
\end{figure}

\vspace{0.25cm}
\noindent \textbf{Grid Structure}\\
In the unsplit approach, we use a Cartesian grid, with point values in each node of the two-dimensional cell and at the center of each edge. Fig. \ref{fig:grid_cells}(a) provides an illustration.

\vspace{0.25cm}
\noindent \textbf{Conservative Update}\\
As we propose a fully discrete method, we do not employ a direct time integration scheme but follow the idea of integrating the Vlasov equation in both space and time. Applying divergence-theorem and introducing the flux functions $h(f) = vf$ and $g(f) = Ef$, we can derive an update formula for the cell averages $\Bar{f}_{i,j}^{n+1}$ in cell $\mathcal{C}_{i,j}$.

\begin{align*}
	\Bar{f}_{i,j}^{n+1} = \Bar{f}_{i,j}^n - \frac{\Delta t}{\Delta x}\left(
	H_{i+\frac{1}{2},j} -H_{i-\frac{1}{2},j}\right) - \frac{\Delta
		t}{\Delta v}\left( G_{i,j+\frac{1}{2}} - G_{i,j-\frac{1}{2}} \right),
\end{align*}
with numerical fluxes 
\begin{align}
	H_{i+\frac{1}{2}, j} \approx \frac{1}{\Delta t \Delta v}\int_{t^n}^{t^{n+1}} \int_{v_{j-\frac{1}{2}}}^{v_{j+\frac{1}{2}}} h(f(x_{i+\frac{1}{2},v},t)) \,dv \,dt,
\end{align}
and
\begin{align}
	G_{i, j+\frac{1}{2}} \approx \frac{1}{\Delta t \Delta x} \int_{t^n}^{t^{n+1}} \int_{x_{i-\frac{1}{2}}}^{x_{i+\frac{1}{2}}} g(f(x, v_{j+\frac{1}{2}}, t)) \,dx \,dt.
\end{align}

These integrals are approximated using Simpson's rule, i.e.
\begin{equation}\label{G_ij}
	\begin{aligned}
		G_{i,j+\frac{1}{2}} = \frac{1}{36}\Big(&E_{i-\frac{1}{2}}^nf_{i-\frac{1}{2}, j + \frac{1}{2}}^{n} + 4E_i^n f_{i, j + \frac{1}{2}}^n + E_{i+\frac{1}{2}}^n f_{i+\frac{1}{2}, j + \frac{1}{2}}^n \\
		+ 4&E_{i-\frac{1}{2}}^{n+\frac{1}{2}}f_{i-\frac{1}{2}, j + \frac{1}{2}}^{n + \frac{1}{2}} + 16E_{i}^{n+\frac{1}{2}}f_{i, j + \frac{1}{2}}^{n + \frac{1}{2}} + 4E_{i+\frac{1}{2}}^{n+\frac{1}{2}}f_{i+\frac{1}{2}, j + \frac{1}{2}}^{n + \frac{1}{2}} \\
		+ &E_{i-\frac{1}{2}}^{n+1}f_{i-\frac{1}{2}, j + \frac{1}{2}}^{n + 1} + 4E_{i}^{n+1}f_{i, j + \frac{1}{2}}^{n + 1} + E_{i+\frac{1}{2}}^{n+1}f_{i+\frac{1}{2}, j + \frac{1}{2}}^{n + 1}\Big),
	\end{aligned}
\end{equation}
and for $H_{i+\frac{1}{2},j}$ respectively.
Subsequently, $f^{n+\frac{1}{2}}, f^{n+1}$, $E^{n+\frac{1}{2}}, E^{n+1}$ are point values of the density distribution function and the electric field at the intermediate time-level $t_{n+\frac{1}{2}}$ and the new time-level $t_{n+1}$.
This approach poses two questions: First, how to evolve the point values, and second, how to evolve the electric field within a time step.

\vspace{0.25cm}
\noindent \textbf{Update of Interface Values}\\
To describe the evolution of point values on the boundaries of each grid cell, we compute the characteristics of the Vlasov equation. In other words, we solve the ODE
\begin{equation}\label{eqn:odeVP}
	\begin{split}
		x'(t) & =  v(t), \\
		v'(t) & =  -E(x(t),t)
	\end{split}
\end{equation}
backwards in time with initial data defined at $t_{n+\frac{1}{2}}$ or $t_{n+1}$. 
The electric field at intermediate times needs to be computed in advance, depending on the stage of the numerical ODE-solver.
This component of the algorithm is addressed in the following subsection.
By solving the ODE backwards in time, we approximately know the trajectory of the point values within a time step along which the value does not change. Subsequently, evaluating these points at time $t_n$ in the respective reconstruction gives approximations of the point values at $t_{n+\frac{1}{2}}$ and $t_{n+1}$ needed in the numerical fluxes above.

\vspace{0.25cm}
\noindent \textbf{Treatment of the Electric Field}\\
In each time step we initially solve the Poisson problem with right-hand-sides
\begin{align*}
	1 - \rho_{i}^n \approx 1 - \int f(x_{i},v, t_n) \,dv,
\end{align*}
using higher-order finite difference schemes with mesh width $\frac{\Delta x}{2}$.
To give a full description of the method, we also have to discuss the approximation of the electric field within a time step $[t_n,t_{n+1}]$ at points $t_{n+\frac{1}{2}}, t_{n+1}$ and every additional point in time arising in the numerical ODE-solver.
In \cite{RS2011}, the authors proposed to use a Taylor expansion in time of the electric field. To achieve third-order accuracy in time, we use
\begin{align}\label{eqn:taylor}
	E(x, t_n + \tau) = E(x, t_n) + \tau \partial_t E(x,t_n) + \frac{1}{2!} \tau^2 \partial_{tt} E(x,t_n) + \mathcal{O}(\tau^3)
\end{align}
for $\tau \in [0, \Delta t]$.
However, direct computations of the derivatives of the electric field are not available, necessitating an alternative approach for their approximation. Moments of the Vlasov equation provide us the following representations:
\begin{align}
	\partial_t E(x,t) &= - \rho u (x,t),\label{eq:E_t} \\
	\partial_{tt} E(x,t) &=  \partial_x \mathbb{E}(x,t_n) - \rho(x,t_n)E(x,t_n), \label{eq:E_tt}
\end{align}
with $\rho u(x,t) := \int v f \, dv$ and $\mathbb{E} := \int v^2 f \,dv$.
The inherent structure of our Cartesian-grid Active Flux method, suggests to approximate $\rho$ and $\rho u$ using Simpson's rule, while the spatial derivative in Eq.~\eqref{eq:E_tt} is carried out by a finite difference scheme on a grid with mesh width $\frac{\Delta x}{2}$.

\subsection{A split-step Active Flux Method for Vlasov--Poisson} 
This subsection briefly summarizes the third-order accurate split-step Active Flux approach presented in \cite{HENSEL2025114294}. 

\vspace{0.25cm}
\noindent \textbf{Grid Structure}\\
The previously discussed two-dimensional Active Flux grid needs to be adjusted to constitute the one-dimensional Active Flux cell and to contain the required degrees of freedom of averages and point values at the interfaces along every direction.
As a consequence of the resulting grid we need to take into account the different types of quantities and provide suitable ways for conversion between them.
The occurring one-dimensional line averages $\Hat{f}_{i+\frac{1}{2}, j}^n, \Hat{f}_{i, j+\frac{1}{2}}^n$ as well as the two-dimensional cell averages $\bar{f}^n_{i,j}$ are initialized from a grid of point values using Simpsons's rule. 

\vspace{0.25cm}
\noindent \textbf{Time Splitting} \\
High-order operator splitting schemes can be constructed on the basis of the fundamental second-order Strang splitting that has originally been introduced as a splitting scheme for the integration of the Vlasov--Poisson system in \cite{cheng1976integration}.  
A standard choice for a fourth-order splitting scheme (here denoted as Yoshida) was introduced in \cite{YOSHIDA1990262} and consists of three consecutive Strang splitting steps (overall 9 steps) with the individual step sizes $\gamma_1 \Delta t$, $\gamma_2 \Delta t$, $\gamma_1 \Delta t$, where $\gamma_1 = 2 / (2-2^{1/3})$ and $\gamma_2 = -2^{1/3} / (2-2^{1/3})$. 
Blanes and Moan \cite{blanes2002practical} introduced a 13-step fourth order operator splitting scheme (here denoted by BM) which is optimized for efficiency. We use both schemes in the numerical experiments to demonstrate the influence on the temporal error in our Active Flux methods.\\

\noindent \textbf{Conservative Update} \\
For the third-order accurate update of the cell average $\bar{f}_{i,j}$ we require computation of the flux integrals Eq.~\eqref{G_ij} with the interface data obtained after the respective directional update steps. 
Then, the conservative update for a time step $\Delta t$ along the $x$-direction is performed as
\begin{align}
	\bar{f}^{*n+1}_{i,j} = \bar{f}^{n}_{i,j} - \frac{\Delta t}{\Delta x} (H_{i+\frac{1}{2},j} - H_{i-\frac{1}{2},j})  \label{eq:FV_xstep}
\end{align}
where we use an asterisk to emphasize the fact that the cell-average has only been updated along one direction and hence represents an intermediate result during the application of the operator splitting scheme. The $v$-step is carried out analogously. 
Due to the fact that the grid for the split-step formulation contains the line-averages over the edges, the point values $f^{n}_{i,j+\frac{1}{2}}$ in the flux integral expression formula must be computed from the one-dimensional reconstruction over the corresponding edges via
\begin{align}
	f^{n}_{i+\frac{1}{2},j} = \frac{1}{4} (6\Hat{f}^{n}_{i+\frac{1}{2},j} - f^{n}_{i+\frac{1}{2},j-\frac{1}{2}} - f^{n}_{i+\frac{1}{2},j+\frac{1}{2}}) \; .
\end{align} 
Note that, due to the splitting, the time evolution of the distribution function $f$ only contains the $x$ or $v$ directional advection in $H_{i+\frac{1}{2},j}$ and $G_{i,j+\frac{1}{2}}$, respectively.   
Furthermore, note that during the one-dimensional directional update in the $v$-direction the electric field is effectively constant as the plasma charge density is an integral over velocity space and hence only a function of $x$. 

\vspace{0.25cm}
\noindent \textbf{Update of Interface Values}\\
On the grid shown in Fig.~\ref{fig:grid_cells}(b) we use consecutive one-dimensional slice-wise updates to evolve the interface values. The edges can thereby be interpreted as one-dimensional Active Flux cells as introduced in \cite{ER2011a}. Update procedures for the DOFs of a one-dimensional Active Flux cell can be compactly expressed as update formulas derived in \cite[p.~ 30,32]{maeng2017advective}.

\vspace{0.25cm}
\noindent \textbf{Treatment of Electric Field}\\
The electric field is obtained by solving the Poisson equation before each $v$-step.  
Integration along the velocity dimension is then carried out using the $\hat{f}_{i-\frac{1}{2}, j}$ as well as the cell averages $\bar{f}_{i,j}$. In order to utilise a nodal Poisson solver, we also need to convert the one-dimensional line averages $\hat{\rho}_i$ into the corresponding point values $\rho_i$ by evaluating the reconstruction at the central location
\begin{align}
	\rho_{i+\frac{1}{2}} = \Delta v \sum_{j=1}^{N_v} \hat{f}_{i+\frac{1}{2},j}, \quad \hat{\rho}_{i} = \Delta v \sum_{j=1}^{N_v} \bar{f}_{i,j}, \quad \rho_i = \frac{1}{4}(6 \hat{\rho}_i - \rho_{i-\frac{1}{2}} - \rho_{i+\frac{1}{2}}). 
\end{align}
One-dimensional averages of the electric field are computed, using Simpson's rule, by $\hat{E}_i = \frac{1}{6} (E_{i-\frac{1}{2}} + 4 E_i + E_{i+\frac{1}{2}})$, which yields the required DOF for the conservation update of the two-dimensional cell averages during the $v$-step of the splitting scheme. 

\subsection{Comparison of the two Methods}
The unsplit Active Flux method uses a uniform grid that includes point values at cell corners and edge midpoints, as well as cell averages. The split-step approach, on the other hand, yields a grid structure combining point values, edge-based line averages, and two-dimensional cell averages.
For the unsplit method, interface point values are evolved in time by solving ordinary differential equations.
The split-step method instead updates interface quantities through successive one-dimensional advection steps with constant velocities within each splitting stage.
Within the unsplit formulation, the time evolution of the electric field is explicitly incorporated into each time step. A single Poisson solve is performed per time step, and the resulting time-dependent field is used in both the ODE solver for point values and the flux computation. The split-step method, in contrast, recomputes the electric field at each velocity update stage, generally requiring multiple Poisson solves per time step. Note that in \cite{HENSEL2025114294}, the convergence error was not dominated by the time-splitting error, and Strang splitting delivered  third-order accuracy overall. Applying second-order Strang splitting reduces the number of Poisson solves required per time step to one.
Numerical fluxes in the unsplit method are obtained from a single two-dimensional second-order polynomial reconstruction and are evaluated once per spatial direction and time step. The split-step approach, on the other hand, relies on multiple one-dimensional reconstructions, with corresponding flux evaluations performed at every stage of the splitting algorithm.
To quantify this, the unsplit method uses a two-dimensional polynomial locally and therefore evaluates $9$ terms per point value. Whereas for the split-step method the one-dimensional reconstruction with $3$ terms is evaluated multiple times successively. Assuming the application of fourth order accurate Yoshida splitting, consisting of three successive Strang splitting steps, we get $3^3 = 27$ term evaluations per time step for each point-value. However, to explore the expandability of the two methods, we can extend this logic to an arbitrary number of dimensions, $d$. One can then approximate the scaling of the term evaluations for a naive expansion of the unsplit method on a tensor product grid as proportional to $3^d$, whereas the split method would require proportional to $3 \cdot N$ term evaluations, where $N$ is the number of directional one-dimensional updates per timestep for a chosen operator splitting scheme in $d$ dimensions. For the Yoshida splitting e.g. the number of update operations in $d$-dimensions scales as $N=3\cdot(2d-1)$, hence for $d\ge 4$, the split method would be more cost-effective than the unsplit method according to this metric. Note that using more accurate splitting schemes such as the optimised Blanes/Moan splitting \cite{blanes2002practical} furthermore increases the costs , as $N$ increases, shifting the raw break-even-point of term evaluations to higher-dimensionalities while of course reducing the splitting induced errors, which also influences the efficiency considerations. Furthermore note that for a third order unsplit method it is not required to use the full tensor-product grid of $3^d$ DOF per cell \cite{arnold2011serendipity}. 
Time integration in the unsplit method is performed in a fully coupled manner, accounting for the simultaneous evolution of all degrees of freedom and the electric field. The split-step method decouples the dynamics through operator splitting, which simplifies the treatment of individual substeps but introduces additional stages within each time step.
For higher-dimensional problems, extending the split-step method is straightforward and has been discussed in \cite{grunwald2025solving}. The underlying framework readily accommodates various splitting algorithms and has been implemented in the \texttt{muphyII} code \cite{allmann2024muphyii}. 
For the unsplit method however, such an extension is not immediate, and it remains unclear how the approach can be generalized.
One major problem would be the representation of the electric field as a function in time.
Although one can formulate a Taylor expansion of the electric field, appropriate expressions of the time-derivatives are not given as Eq.~\eqref{eq:E_t}, \eqref{eq:E_tt} hold only in the one-dimensional case.

\medskip

\section{Numerical Results}\label{sec:3}
In this subsection, we present a numerical comparison of the previously described methods. We perform this study on two standard computational plasma physics problems: the two-stream instability and strong Landau damping. For both Active Flux methods we use the same \textit{CFL} condition $\frac{\Delta t v_{max}}{\Delta x} = \frac{2}{\pi} \approx 0.63661977$ to determine the time step. 
Additionally, we use a split-step \textit{PFC} method \cite{FSB2001} with high grid resolution ($N_x = N_v = N_{ref} = 512$) combined with the fourth order splitting scheme by Blanes/Moan and a small time step (by setting the \textit{CFL} here to $\frac{1}{5\pi} \approx 0.06366197$) as a reference to verify our numerical results.

\subsection{Two Stream Instability}
We consider the following initial electron distribution function
\begin{align}
	f_e(x, v, t=0) = \frac{1}{2  \sqrt{2\pi}} \left[ \exp\left(-\frac{(v-v_0)^2}{2} \right) + \exp\left(-\frac{(v+v_0)^2}{2}\right)  \right] (1 + A \cos(kx))
\end{align}
on the periodic phase space volume $[-5\pi,5\pi]_x \times [-10, 10]_v$. The beam velocities are set to $v_0 = 3$ and the perturbation parameters to $k=0.2$ and $A=10^{-3}$. 
For this example we perform a long-time accuracy study for the distribution function simulated by both previously described Active Flux methods. In order to compare the results to the high-resolution reference, we scale the \textit{PFC} solution down to a desired lower resolution by averaging over neighboring cells. 

\begin{figure}
	\centering
	\includegraphics[width=\linewidth]{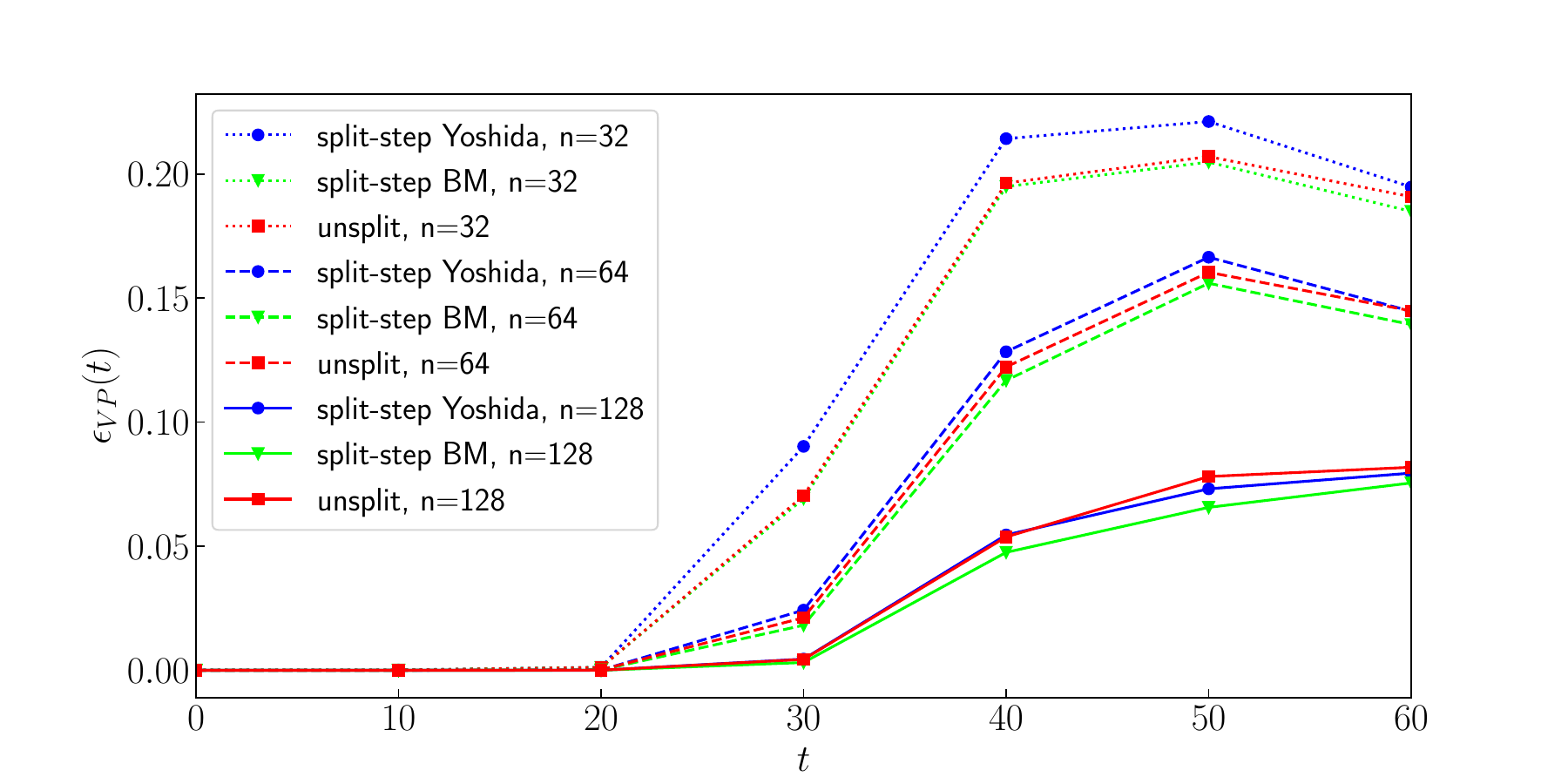} 
	\caption{Error $\epsilon_{VP}$ over time for the two stream instability with varying grid resolutions $(N_x,N_v)=(n,n)$.}
	\label{fig:TS_error_over_time}
\end{figure}

To quantify the error of our numerical results we consider the following relative-error-measure:
\begin{align}
	\epsilon_{VP} = \frac{\sum_{i=1}^{N_x} \sum_{j=1}^{N_v} |\Tilde{\bar{f}}_{i,j}- \bar{f}_{i,j}|}{\sum_{i=1}^{N_x} \sum_{j=1}^{N_v} |\Tilde{\bar{f}}_{i,j}|}.
\end{align}
where $\Tilde{\bar{f}}_{i,j}$ is the scaled down reference electron distribution and $\bar{f}_{i,j}$ the numerical result obtained by one of the Active Flux methods.
Note that a comparison using a similar methodology, with a reference solution obtained by simulating with a small CFL, has been performed between split and unsplit SLDG schemes on the Vlasov--Poisson system in \cite{cai2022comparison}.
In Fig.~\ref{fig:TS_error_over_time}, we present this relative error on grids with $N_x = N_v = 32, 64, 128$ cells. 
During the initial phase of the instability, up until $t\approx20$, smaller errors are observed due to the smoothness of the initial conditions. As the instability develops, the solution structure grows in complexity. This results in constant magnification of the numerical errors up until $t\approx40$.  After $t\approx40$, the grid resolution can no longer resolve the finer structures of the solution, and numerical dissipation becomes the dominant effect. 

Both the split and unsplit methods show the described error behavior, as well as similar scaling of error reduction with increasing grid resolution. We find that the applied splitting scheme visibly influences the error for the split-step method showing that the temporal error dominates. Although the fourth-order Yoshida splitting produces larger errors than the unsplit scheme, the optimised splitting by Blanes/Moan offers slightly better results than the unsplit scheme.

\begin{figure}
	\centering
	\includegraphics[width=1\textwidth]{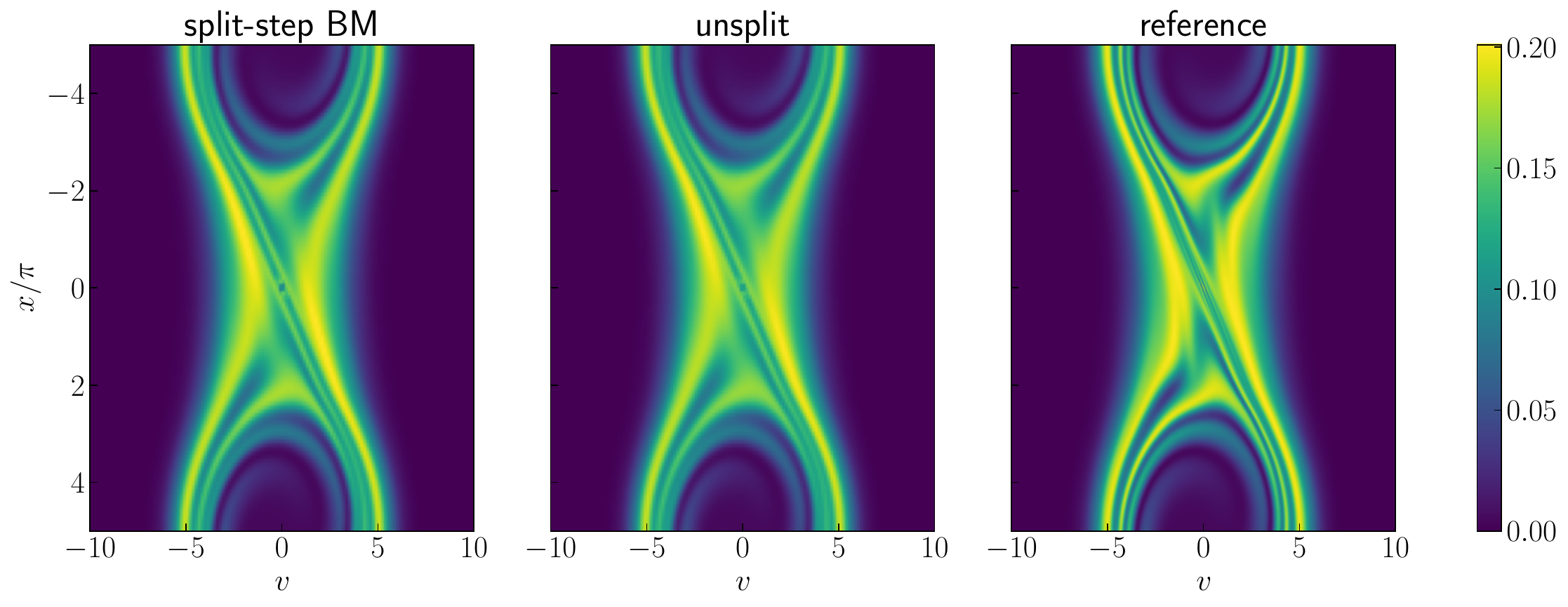}
	\caption{Snapshots of the plasma distribution functions with $(N_x, N_v)=(128, 128)$ compared to a high resolution reference at time $t= 50$.}
	\label{fig:snapshots_TS}
\end{figure}

Additionally, Fig.~\ref{fig:snapshots_TS} shows snapshots of the electron distribution function as simulation results of the two Active Flux methods, compared to the high-resolution \textit{PFC} result at $t = 50$. Note that the reference solution is computed on a grid with $512\times 512$ cells, while both Active Flux methods use a grid with $128 \times 128$ cells. In the absence of an analytical solution, the reference snapshot serves as a benchmark for displaying the correct solution structure.

\subsection{Strong Landau Damping}
\begin{figure}
	\centering
	\begin{subfigure}[b]{0.49\textwidth}
		\centering
		\includegraphics[width=\textwidth]{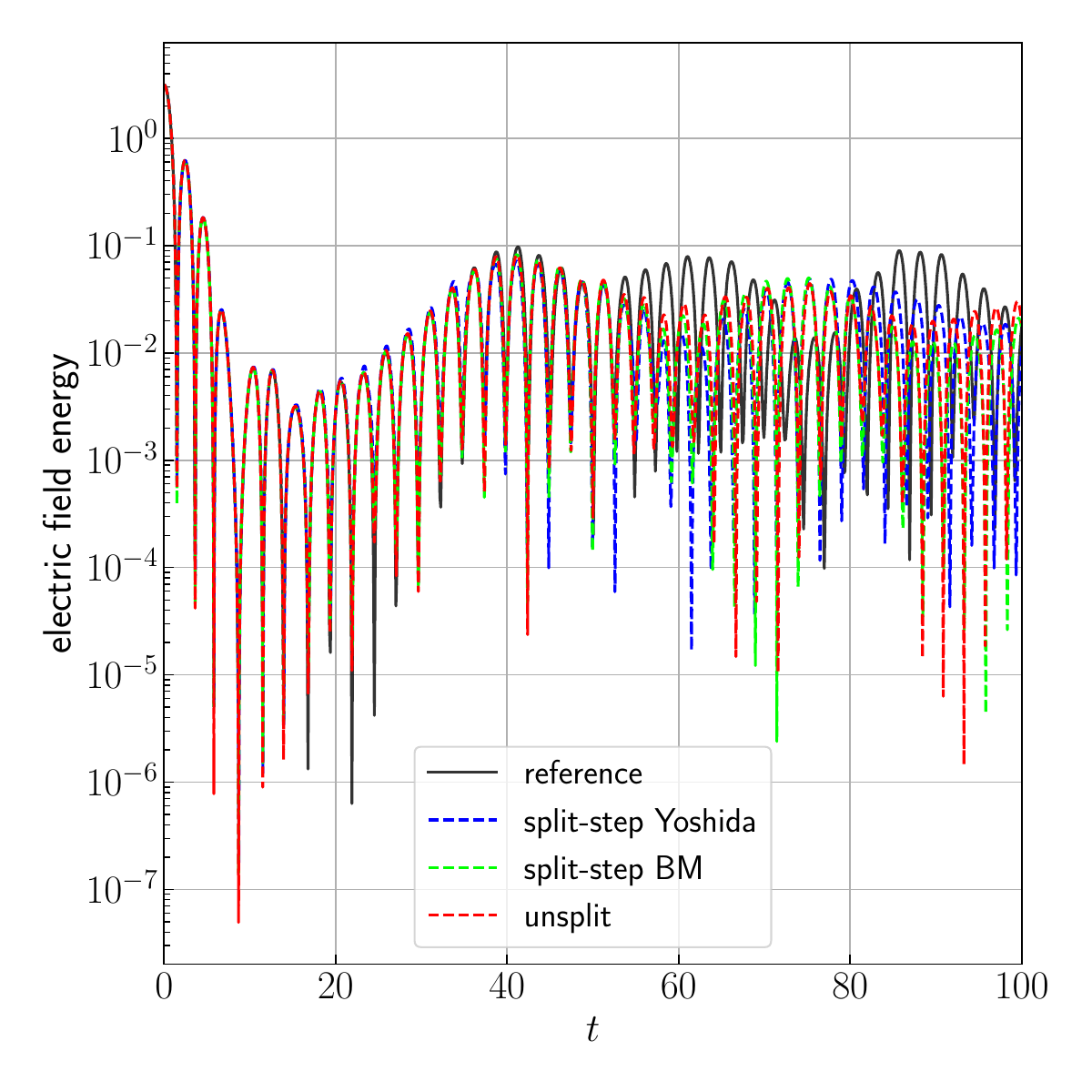}
		\caption{$(N_x, N_v) = (64, 64)$}
		\label{fig:SLD_a}
	\end{subfigure}
	\hfill
	\begin{subfigure}[b]{0.49\textwidth}
		\centering
		\includegraphics[width=\textwidth]{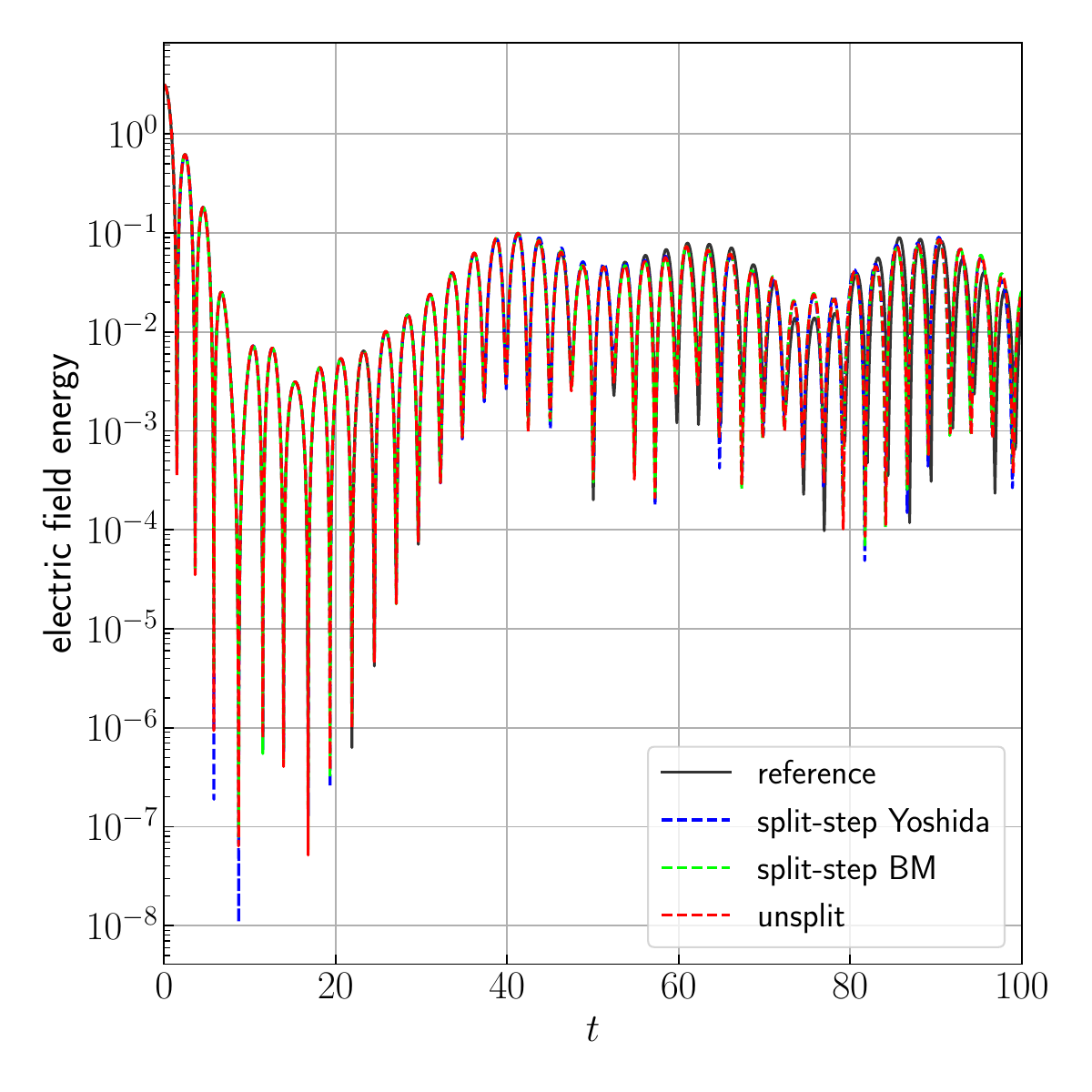}
		\caption{$(N_x, N_v) = (128, 128)$}
		\label{fig:SLD_b}
	\end{subfigure}
	\caption{Long-time behavior of the electric field energy compared to a high-resolution reference result}
	\label{fig:strong_LD_electric_field}
\end{figure}

In order to examine the long-time behavior of the electric field for both methods we furthermore consider the example of strong Landau damping, with the electron distribution function 
\begin{align}
	f_e(x, v, t=0) = \frac{1}{\sqrt{2 \pi}} \exp{\left(-\frac{v^2}{2} \right)}  (1 + A \cos(kx)),
	\label{eq:LD_initial_condition}
\end{align}
initialised on the phase space domain $[ -2\pi, 2\pi]_x \times  [-5, 5]_v$, with periodic boundaries in both direction.
We set the parameters $A = 0.5$ and $k=0.5$.
In Fig.~\ref{fig:strong_LD_electric_field} we present the evolution of the discrete electric field energy $\frac{1}{2} \Delta x \sum_{i=1}^{N_x} \hat{E}_i^2$ for both Active Flux methods and the high-resolution reference solution. 
We use a highly resolved \textit{PFC} solution as reference of the damping behavior. 
As in the previous comparison, the results for both Active Flux methods are in close agreement.
On the coarser grid, we see that both Active Flux methods behave accordingly for the first half of the simulation (Fig.~\ref{fig:SLD_a}). After $t \approx 50$ the Active Flux methods do not match the the electric energy of the reference solution.
With increased resolution (Fig.~\ref{fig:SLD_b}), both Active Flux methods show an improved long-term behavior that is in accordance with the reference solution.  

\section{Conclusion}
In this paper the two existing Active Flux based approaches for the kinetic Vlasov--Poisson system have been compared algorithmically and quantitatively in terms of their numerical results on two benchmark problems from electrostatic plasma physics. 
While the unsplit method is fully discrete and delivers slightly better numerical results than the split-step approach in our comparisons to the reference on coarser grids, it is also cheaper for the two-dimensional case considered here. However, the split-step approach has slightly better results on finer grids and can be improved by using more accurate splitting schemes that come with further increase of the computational cost and more easily expanded to higher dimensions and the electromagnetic Vlasov--Maxwell system.

\section*{Acknowledgements}
We acknowledge funding from the German Science Foundation DFG through the research unit ``SNuBIC'' DFG-FOR5409, project id 463312734 and project id 530709913.

\bigskip\bigskip

\printbibliography

\end{document}